\documentstyle[12pt, twoside]{article}
\input amssym.def
\input amssym.tex 

\newenvironment{demo}[1]%
{\vskip-\lastskip\medskip
  \noindent
  {\em #1.}\enspace
  }%
{\qed\par\medskip
  }

\newcommand{\qed}{
  \strut\hfill
  \mbox{$\Box$}
  }

{\vskip-\lastskip
  \begin{trivlist}
    \item[]
      {\bf Axiom #1 }
      }%
{\end{trivlist}
  }

\newtheorem{theorem}{Theorem}[section]

\newtheorem{corollary}{Corollary}[section]

\newtheorem{lemma}{Lemma}[section]

\newtheorem{remark}{Remark}[section]

\newtheorem{definition}{Definition}[section]

\newtheorem{proposition}{Proposition}[section]

\newcommand{\bz}{
  \sum_{n \in \Z} b(n) z^{-n }
  }

\newcommand{\cz}{
  \sum_{n \in \Z} c(n) z^{-n -1}
  }

\newcommand{\C}{
 \Bbb C
  }

\newcommand{\D}{
  \cal D
  }

\newcommand{\End}{
  \mbox{End }
  }

\newcommand{\F}{
  {\cal F}
  } 

\newcommand{\Fl}{
  {\cal F}^l
  }

\newcommand{\Fo}{
  {\cal F}^0
  }

\newcommand{\Fbaro}{
  {\overline{\cal F} }^0
  }

\newcommand{\Halpha}{
  {\cal H}_{\alpha}
  }

\newcommand{\Hl}{
  {\cal H}_l
  }
  
\newcommand{\Ho}{
  {\cal H}_0
  }
  
\newcommand{\hD}{
  { \widehat{\cal D} }
  }

\newcommand{\hf}{
  \frac12
  }

\newcommand{\Ltw}{
   {\cal L}(t, w)
  }

\newcommand{\Moo}{
   {\cal M} (0, 0)
  }

\newcommand{\Mtw}{
   {\cal M}(t, w)
  }

\newcommand{\module}{
  M = \bigoplus_{n \in \Z_{+}} M_n
  }

\newcommand{\ps}{
  \psi (z)
  }

\newcommand{\strz}{Res_z
  \left(
    Y(a,z) \frac{(z+1)^{\wt a}}{z} b
  \right)
  }

\newcommand{\strrz}{Res_z
  \left(
    Y(a,z) \frac{(z+1)^{\wt a}}{z^2} b
  \right)
  }

\newcommand{\SUM}{
  \sum_{n\in \Bbb Z}
  }

\newcommand{\Tbc}{
  :\partial b(z) c(z):
  }

\newcommand{\UW}{
  \/{\cal U}({\cal W}_3)
  }

\newcommand{\V}{
  V = \bigoplus_{n \in \Z_{+}} V_n
  }

\newcommand{\vac}{
  | 0 \rangle
  }

\newcommand{\vacbc}{
  | bc \rangle
  }

\newcommand{\vactw}{
  | t, w \rangle
  }

\newcommand{\vacuum}{
  | \alpha \rangle
  }

\newcommand{\vacoo}{
  | 0,0 \rangle
  }

\newcommand{\VM}{
  {\cal VW}_{3,c}
  }

\newcommand{\VMtwo}{
  {\cal VW}_{3,-2}
  }

\newcommand{\W}{
  { \cal W }_{1+\infty}
  }

\newcommand{\Wbc}{
  \frac{1}{\sqrt{6}}
   \left( :\partial^2 b(z) c(z): - :\partial b(z) \partial c(z): \right)
  }

\newcommand{\Wc}{
  {\cal W}_{{1 + \infty}, c}
  }

\newcommand{\wt}{
  \mbox{wt }
  }

\newcommand{\Wth}{
  {\cal W}_3
  }

\newcommand{\Wtho}{
  {\cal W}_{3,0}
  }

\newcommand{\Wthp}{
  {\cal W}_{3, +}
  }

\newcommand{\Wthpm}{
  {\cal W}_{3, \pm}
  }

\newcommand{\Wthc}{
  {\cal W}_{3, c}
  }

\newcommand{\Wthtwo}{
  {\cal W}_{3, -2}
  }

\newcommand{\WGN}{
  {\cal{W} } (gl_N)
  }

\newcommand{\Winfone}{
  {\cal W}_{{1 + \infty}, -1}
  }

\newcommand{\WN}{
  {\cal W}_N
  }

\newcommand{\Z}{
  \Bbb Z
  }

\newcommand{\zwone}{
  (z-w)
  }

\newcommand{\zwtwo}{
  (z-w)^2
  }

\newcommand{\zwthree}{
  (z-w)^3
  }

\newcommand{\zwminusone}{
  \frac{1}{z-w}
  }

\newcommand{\zwminustwo}{
  \frac{1}{(z-w)^2}
  }

\begin{document}
\title{ 
 Classification of irreducible
 modules \\ of $\Wth$ algebra
  with $c = -2$
  }
\author{
  %
  Weiqiang Wang\\
\\{\small Max-Planck Institut fur Mathematik, 53225 Bonn, Germany}
\\{\small E-mail: wqwang@mpim-bonn.mpg.de}}

\date{}

\maketitle
\begin{abstract}
  We construct irreducible modules $V_{\alpha}, \alpha \in \C$
  over $\Wth$ algebra with $c = -2$ in terms of
  a free bosonic field. We prove that these modules
  exhaust all the irreducible modules of
  $\Wth$ algebra with $c = -2$. 
  Highest weights of modules $V_{\alpha}, \alpha \in \C$ with respect to
  the full (two-dimensional) Cartan subalgebra of
  $\Wth$ algebra are 
  $\left( \hf \alpha(\alpha-1),
  \frac{1}{6} \alpha(\alpha-1)(2\alpha-1) \right)$.
  They are parametrized by points $(t, w)$ on a rational curve
  $ w^2 - \frac{1}{9} t^2 (8t + 1) = 0.$ Irreducible modules of
  vertex algebra $\W$ with $c = -1$ are also classified. 
\end{abstract}
\setcounter{section}{-1}
\section{Introduction}

In the study of two-dimensional conformal field theories
extensions of conformal symmetry play an important role.
The algebraic structures underlying the extended conformal
symmetry are usually known as $\cal W$-algebras in 
literatures (see \cite{BS, FF} and references 
therein). Mathematically $\cal W$-algebras can be put into
the general framework of the theory of vertex algebras 
formulated first by Borcherds, 
cf. e.g. \cite{B1, FLM, DL, LZ, FKRW, K, B2}. 

In contrast to vertex algebras associated to 
the Virasoro algebra, $\cal W$-algebras
such as $\WN$ algebras, have the feature that
non-linearity terms appear in the operator product
expansion of two generating fields, namely 
the commutator of two generators contains
non-linear terms expressed by these generators themselves.
Mainly due to the non-linear nature of $\cal W$-algebras,
the study of their representation theory 
has been difficult and very non-trivial. Even
the understanding of representation theory
of the Zamolodchikov $\Wth$ algebra \cite{Za}, 
which is the simplest example of
$\cal W$-algebras beyond the Virasoro algebra, 
is far from satisfactory (see however \cite{BMP}). 

Apart from the Virasoro generators $L_n$, $n \in \Z$,
$\Wth$ algebra has an additional set of generators
$W_n$, $n \in \Z$. Denote by $\UW$
the corresponding universal enveloping algebra.
Define two generating series
$T(z) = \SUM L_n z^{-n-2}$ and 
$W(z) = \SUM W_n z^{-n-3}$. 

It is well known that
the vacuum module $\VM$ with central charge $c$
carries a vertex algebra structure. 
For a generic central charge $c$, $\VM$ is an irreducible 
representation of $\UW$. In this case, the representation
theory of vertex algebra $\VM$ is the same as 
that of $\UW$. For a non-generic central charge~$c$,
$\VM$ is reducible and admits
a unique maximal proper $\UW$-submodule $ I$ and thus
a unique irreducible quotient, which is denoted 
by $\Wthc$. $\Wthc$ inherits a vertex algebra structure from 
$\VM$. Representation theory of $\Wthc$ with
non-generic central charge $c$ becomes highly non-trivial
since a module $M$ of $\UW$ can be regarded as a module
of $\Wthc$ if and only if the Fourier components of
any field corresponding to any vector in $ I$
annihilates the whole $M\/$. 

In \cite{W2}, in studying the vertex algebra $\W$ with central 
charge $-1$ (denoted by $\Winfone$) 
we explicitly constructed a number of
irreducible modules 
of $\Wthtwo$ parametrized by integers
and obtained full character formulas for these
modules. We showed that 
the vertex algebra $\Winfone$ is isomorphic 
to a tensor product of $\Wthtwo$ and a Heisenberg vertex
algebra generated by a free bosonic field
by using Friedan-Martinec-Shenker bosonization technique \cite{FMS}.

In this paper, we will continue the study of
representation theory of $\Wthtwo$ and
$\Winfone$. Note that $-2$ is
a non-generic central charge for $\Wth$ algebra.
We will explicitly construct irreducible 
modules $V_{\alpha}, \alpha \in \Bbb C$ of $\Wthtwo$ in terms of
a free bosonic field. Then by locating 
key singular vectors in $\VMtwo$ and then applying
Zhu's machinery \cite{Z} to our case we
are able to prove that $V_{\alpha}, \alpha \in \Bbb C$ exhaust
all the irreducible modules of $\Wthtwo$. It 
turns out that the set of all irreducibles of
$\Wthtwo$ has an elegant description:
highest weights of these irreducible modules
are parametrized by points of a rational curve 
defined by $w^2 - \frac{1}{9} t^2 (8t + 1) = 0.$
Combining with our
results in \cite{W2} we also construct and classify all
the irreducible modules of $\Winfone$.
This latter classification result disproves 
a conjecture of Kac and Radul \cite{KR2}.

Let us explain in more detail.
Given a pair of $bc$ fields 
$b(z) = \bz$ and $c(z) = \cz$, we construct
a Fock space $\F$ generated by the vacuum
vector $ \vacbc$, satisfying
$$ b( n+1 ) \vacbc = 0, 
   \quad  c(n) \vacbc = 0, \quad n \geq 0. $$
Then $j (z) = :b(z) c(z):= \SUM j_n z^{-n-1}$ 
is a free bosonic field. Take a scalar field $\psi (z)$
such that $j(z) = \partial \psi (z)$.
Denote by $\Halpha$ the
Fock space of the Heisenberg algebra $\{ j_n, n \in \Z \}$ 
with vacuum vector $\vacuum$ satisfying
$$j_n \vacuum = \alpha \delta_{n, 0} \vacuum, \quad n \geq 0.$$
It is observed in \cite{BCMN, W2} that
the fields 
$$ T(z) = \Tbc, \quad W(z) = \Wbc , $$
satisfy the $\Wth$ operator product expansions
with central charge $-2$. We can rewrite the fields
$T(z)$ and $W(z)$ in terms of $j(z)$ by means of
boson-fermion correspondence (Proposition \ref{prop_bosonw3}). 
$\Fo$ is thus isomorphic to $\Ho$.
It is shown \cite{W2}
that the simple vertex algebra $\Wthtwo$ is a vertex subalgebra
of $\Fo$ and can be identified explicitly inside $\Fo$. 
Denote by $V_{\alpha}$ the irreducible quotient
of the $\Wthtwo$-submodule of $\Halpha$ generated 
by the highest weight vector $\vacuum$ in $\Halpha$.
Let $\widetilde{W} (z) \equiv \SUM \widetilde{W}_n z^{-n-3}
= \hf \sqrt{6} W(z)$.
We will show that the highest weight of $V_{\alpha}$ 
with respect to the full Cartan subalgebra 
$\{ L_0, \widetilde{W}_0 \}$ of $\Wthtwo$ is
\begin{eqnarray}
  \left( \hf \alpha(\alpha-1), \frac{1}{6} \alpha(\alpha-1)(2\alpha-1) \right).
 \label{eq_wt}
\end{eqnarray}

To show that the above irreducible modules 
$V_{\alpha}, \alpha \in \C$ exhaust all the irreducible
modules of $\Wthtwo$, we invoke a powerful machinery 
due to Zhu in the general theory of vertex algebras \cite{Z}.
Zhu constructed an associative algebra $A(V)$ for 
any vertex algebra $V$ such that irreducible 
modules of the vertex algebra $V$
one-to-one correspond
to irreducible modules  
of the associative algebra $A(V)$. Zhu's constructions
were generalized to vertex superalgebras in \cite{KW}.
By construction, the Zhu associative algebra $A(V)$ is
a certain quotient of $V$. We denote by $[a]$
the image in $A(V)$ of $a \in V$.
By studying the associative algebra
$A(V)$, one can often obtain useful information on 
highest weights of modules over $V$. 

We show that the Zhu associative
algebra $A(\VMtwo )$ is isomorphic to a polynomial
algebra ${\Bbb C} [t, w]$, where $t$ and $w$ 
correspond to elements $[ L_{-2} \vac ]$ and
$[ \widetilde{W}_{-3} \vac ]$ in $A(\VMtwo )$
respectively. Using some explicit results on
singular vectors in $\VMtwo$, we further show that
Zhu associative algebra
$A(\Wthtwo )$ is isomorphic to
(some quotient of) the quotient algebra ${\C} [t, w]/<f(t, w)>,$
where $\/<f(t, w)>\/$ is the ideal of ${\C} [t, w]$
generated by the polynomial $f(t, w) = w^2 - \frac{1}{9} t^2 (8t + 1).$
This means that a necessary condition
for any irreducible module
of the vertex algebra $\VMtwo$ to be
a module of $\Wthtwo$ is that
its highest weight $(t, w)$ has to satisfy
the equation 
\begin{equation}
  w^2 - \frac{1}{9} t^2 (8t + 1) = 0.
 \label{eq_elliptic}
\end{equation}

We observe that all the solutions to the 
equation above can be written as of the form
(\ref{eq_wt}). But we have already constructed  
irreducible modules $V_{\alpha}$ of $\Wthtwo$ with
a highest weight of any such form. This shows that
$V_{\alpha}$ $ (\alpha \in \C)$
are all irreducible $\Wthtwo$-modules
and their highest weights are parametrized by
points on the rational curve defined by
the equation (\ref{eq_elliptic}). The equation (\ref{eq_elliptic})
as a necessary constraint
on the highest weights of irreducible $\Wthtwo$-modules
was anticipated in \cite{H, EFHHNV} by some other arguments
\footnote{I thank A. Honecker for mentioning these references to me.}.

In the remaining part of this paper we classify
all the irreducible modules of
vertex algebra $\Winfone$. Recall in
our paper \cite{W2}, we have shown that
the vertex algebra $\Winfone$ is isomorphic to
a tensor product of the vertex algebra $\Wthtwo$
and a Heisenberg vertex algebra generated by
a free bosonic field. Therefore the
classification of irreducible modules
$\Winfone$ follows from our classification of
irreducible modules of $\Wthtwo$
and the well-known description of all irreducible 
modules of a Heisenberg vertex algebra.

This paper is organized as follows. In Section \ref{sect_general},
we recall the definition of a vertex algebra and 
review Zhu's associative algebra 
theory. In Section \ref{sect_wthree}, 
we recall the $\Wth$ algebra and
study the case with central charge $-2$ in some detail.
In Section \ref{sect_realize} we construct 
irreducible modules $V_{\alpha}$ $(\alpha \in \C)$ of $\Wthtwo$
and determine their highest weights.
In section \ref{sect_complete} we calculate
Zhu algebra $A(\Wthtwo)$ and show that the list of irreducible 
modules constructed in Section \ref{sect_realize} 
is complete. In Section \ref{sect_infty} we classify
all irreducible modules of $\Winfone$.

\section{Vertex algebras and Zhu's associative algebra theory}
  \label{sect_general}

Our definition of vertex algebras basically follows \cite{FKRW, K}. 
It is known that
our definition is essentially equivalent to other formulations
in \cite{B1, FLM, LZ, B2}. Though it is not essential to
have a gradation in the definition of vertex algebras,
we choose to keep it
in this paper in order to present Zhu's associative algebra theory.

\begin{definition}
  A vertex algebra consists of the following data:
a $\Z_{+}$-graded vector space $\V$; 
a vector $| 0 \rangle \in V\/$ (called the vacuum vector);
an operator $L_0$ (called the degree operator)
and an operator $T \in \End V $ (called the translation 
operator);
a linear map from $V\/$ to the space of fields 
$a \mapsto Y (a, z) = \sum_{n \in \Z} a_{(n)} z^{- n - 1}
\in \End V [[z, z^{-1}]]$ (called the state-field
correspondence). These data satisfy the following axioms:
 \begin{enumerate}
 \item[(V)] $Y (|0\rangle, z) = I_V, 
   \quad Y (a, z) \vac |_{z = 0} = a$;

 \item[(G)] $L_0 |_{V_n} = n I_{V_n}, \quad
  [ L_0, Y(a, z)] = \partial_z Y(a,z) + Y(L_0 a, z);$

 \item[(T)] $[T, Y (a, z)] = \partial_z Y (a, z), \quad
   T \vac = 0$;

 \item[(L)] $(z - w)^N [Y (a, z), Y (b, w)] = 0 \quad \mbox{for }
   \quad N \gg 0$.
 \end{enumerate}
\end{definition}

For $a \in V_n$, $n$ is called the
weight of $a$, denoted by $\wt a$.
Denote by $ o(a) = a (\wt a - 1)$ for homogeneous $ a \in V\/$
and extends by linearity to the whole $V$.
The results of the remaining part of this section are due to
Zhu \cite{Z}. We refer the readers to \cite{Z} for more detail.

\begin{definition}
  Define two bilinear operations $*$ and $\circ $ on $V\/$
 as follows. For $a $ homogeneous, let
  \begin{eqnarray*}
    a*b       & = & \strz,        \\
    a \circ b & = & \strrz,
 \end{eqnarray*}
  then extend to $V\/ \times V\/$ by bilinearity.
  Denote by $O(V)\/$ the subspace of $V\/$
  spanned by elements $ a \circ b $, and by $A(V)$
  the quotient space $V/O(V)$.
 \label{def_mult}
\end{definition}

It is convenient to introduce an equivalence relation $\sim\/$ as
in \cite{W1}.  For $a, b \in V\/$, $a \sim b\/$ means $a - b \equiv 0
\bmod O(V)$. For $f, g \in End\ V\/$, $f \sim g\/$ means $f \cdot c
\sim g \cdot c\/$ for any $c \in V\/$.  Denote by $[a]\/$ the
image of $a\/$ in $V\/$ under the projection of $V\/$ onto $A(V)\/$.

\begin{lemma} 
\begin{enumerate}
\item[1)]  $T + L_0 \sim 0$.
 
\item[2)] For every homogeneous element $a\in V$, and $m\geq n
\geq 0$, one has
\begin{eqnarray*}
  Res_z \left(Y(a,z)\frac {(z+1)^{\wt a+n}}{z^{2+m}}\right) &
  \sim & 0.
\end{eqnarray*}

\item[3)] For homogeneous elements $a,b \in V$, one
  has
  $$ a*b  \sim Res_z
  \left(
    Y(b,z) \frac{(z+1)^{\wt b-1}}{z} a
  \right). $$
\end{enumerate}
\label{lem_keylemma}
\end{lemma} 

\begin{theorem} 
  \begin{enumerate}
  \item[1)] $O(V)$ is a two-sided ideal of $V$ under the
    multiplication $*$. Moreover, the quotient algebra $(A(V),*)$
    is associative.
  \item[2)]  $[1]$ is the unit element of the algebra $A(V)$.
  \end{enumerate}
  \label{theorem_assoc}
\end{theorem} 

In the case that the vertex algebra $V$ 
contains a Virasoro element~$\omega$, i.e. the corresponding
field $\/Y(\omega, z)$ is an energy-momentum tensor field,
we have

\begin{lemma}
  $[\omega] $ is in the center of the associative algebra $A(V)$.
  \label{lem_center}
\end{lemma}

The following proposition follows from the definition of $A(V)$.

\begin{proposition}  
  Let $I$ be an ideal of $V$. Then the associative
  algebra $ A(V/I)$ is isomorphic to $A(V) / [I]$, where $[I]$ is
  the image of $I$ in $A(V)$.
\label{prop_associativity}
\end{proposition}
 
\begin{theorem}  
  \begin{enumerate}
  \item[1)] If $\module$ is a module of the vertex algebra 
  $V$, then the top level 
   $M_0$ of $M$ is a module of the associative algebra
  $A(V)$, with action given as follows: 
  for $[a] \in A(V)$, which is the image of $ a \in V$,
  $[a]$ acts on $M_0$ as $o(a)$.

  \item[2)] Irreducible modules of the vertex algebra $V$ one-to-one
  correspond to irreducible modules of the associative
  algebra $A(V)$ as in 1).
  \end{enumerate}
 \label{th_zhucorresp}
\end{theorem}

We call $A(V)$ the {\em Zhu (associative) algebra} of a vertex
algebra $V$.

\section{$\Wth$ algebra with central charge $-2$}
  \label{sect_wthree}

Denote by $\/{\cal U} (\Wth )$ the quotient of the
free associative algebra
generated by $L_m, W_m,\; m \in \Z$, by the ideal 
generated by the following commutation relations (cf. e.g. \cite{BMP}):
\begin{eqnarray}
{[ L_m, L_n ]} & = &  
    (m-n) L_{m+n} + \frac{c}{12}(m^3 - m) \delta_{m, -n}, \nonumber \\
{[ L_m, W_n ]} & = & 
     (2m - n) W_{m + n},        \label{eq_alg}            \\
{[ W_m, W_n ]} & = &
    (m-n)\Bigl(
           \frac{1}{15}(m+n+3)(m+n+2)  \nonumber     \\
     && \quad\quad\quad\quad\quad\quad\quad\quad\quad - \frac{1}{6} (m+2)(n+2)
         \Bigl) L_{m+n}                 \nonumber    \\
     && + \beta (m-n) \Lambda_{m+n} 
        + \frac{c}{360}m(m^2 -1)(m^2 -4) \delta_{m,-n},   \nonumber    
\end{eqnarray}
where $c \in \Bbb C$ is the central charge, $\beta = 16/(22+5c) $ and
\begin{eqnarray*}
   \Lambda_m = \sum_{n \leq -2} L_n L_{m-n} 
               + \sum_{n > -2} L_{m-n} L_n 
               - \frac{3}{10} (m+2)(m+3) L_m.
\end{eqnarray*}

Denote
\begin{eqnarray*}
   \Wthpm = \{ L_n, W_n, \; \pm n \geq 0 \}, \quad
   \Wtho = \{ L_0, W_0 \}.
\end{eqnarray*}

A Verma module $\/ {\cal M}_c (t, w)$ (or $\/ \Mtw$ 
whenever there is no confusion of central charge)
of $\/{\cal U}(\Wth)\/$ is the
induced module
\begin{eqnarray*}
  \Mtw = 
   {\cal U}(\Wth) \bigotimes_{ {\cal U}( \Wthp \oplus \Wtho) }
    {\Bbb C}_{t, w}
\end{eqnarray*}
where $\/{\Bbb C}_{t, w}$ is the 1-dimensional module of
$\/{\cal U}( \Wthp \oplus \Wtho)$ generated by
a vector $|t,w\rangle$ such that 
\begin{equation}
  \Wthp \vactw = 0, \,L_0 \vactw = t \vactw, \,
   W_0 \vactw = w \vactw.
\end{equation}
$\Mtw $ has a unique irreducible quotient which is denoted by
$\Ltw$ (or ${\cal L}_c(t, w)$ when it is necessary to specify
the central charge). A singular vector in a $\/ {\cal U}(\Wth)$-module
means a vector killed by $\Wthp$. For simplicity, we 
denote the vacuum vector $\vacoo$ by $\vac$ in the case $t = w =0$.
It is easy to see that
$L_{-1} \vac, W_{-1} \vac$, and $W_{-2} \vac$
are singular vectors in $\Moo$. We denote
by $\VM$ the {\em vacuum module} which is by definition the quotient
of Verma module $\Moo\/$ by the $\/ {\cal U}(\Wth)$-submodule 
generated by the singular vectors 
$L_{-1} \vac, W_{-1} \vac$, and $W_{-2} \vac$.
We also call ${\cal L}(0, 0)$ the {\em irreducible vacuum module}.
Let $I$ be the maximal proper submodule of 
the Verma vacuum module $\VM$.
Clearly ${\cal L}(0, 0)$ is the irreducible quotient
of $\VM$. It is easy to see that $\VM$ has a linear basis 
\begin{eqnarray}
  L_{-i_1 -2} \cdots L_{-i_m -2} W_{-j_1 -3} \cdots W_{-j_n -3}\vac,  
       \label{eq_PBW}                   \\
    \quad 0 \leq i_1 \leq \cdots \leq i_m, 
   \quad 0 \leq j_1 \leq \cdots \leq j_n,\quad m,n \geq 0.  \nonumber
\end{eqnarray}

The action of $L_0$ on $\VM$ gives
rise to a principal gradation on $\VM$:
$ \VM = \bigoplus_{n \in \Z} (\VM)_n. $
Introduce the following fields
\begin{equation}
  T(z) = \sum_{n \in \Z} L_n z^{-n-2},  \quad
  W(z) = \sum_{n \in \Z} W_n z^{-n-3}.
\end{equation}
It is well known that the vacuum module
$\VM$ (resp. irreducible vacuum module ${\cal L}(0, 0)$)
carries a vertex algebra structure 
with generating fields $T(z)$ and $W(z)$.
The $\Wth$ algebra with central charge $-2$ we have
been referring to
is the vertex algebra ${\cal L}_{-2}(0, 0)$, which we denote
by $\Wthtwo$ in this paper. 
Fields $T(z)$ and $W(z)$ correspond to the 
vectors $L_{-2} \vac$ and $W_{-3} \vac$ 
respectively. The field corresponding to the vector
$L_{-i_1 -2} \cdots L_{-i_m -2} W_{-j_1 -3} \cdots W_{-j_n -3} \vac$
is 
$$ \/\partial^{(i_1)} T(z) \cdots \partial^{(i_m)} T(z) 
  \/\partial^{(j_1)} W(z) \cdots \partial^{(j_n)} W(z), $$
where $\/\partial^{(i) }\/$ denotes 
$ \/\frac{1}{i!}\partial^i_z $. 

From now on we concentrate on the case of $\Wth$ algebra
with central charge $c = -2$. We can rewrite (\ref{eq_alg})
as the following OPEs in our central charge $ -2$ case:
  \begin{eqnarray}
    T(z) T(w) & \sim & \frac{-1}{(z-w)^4} 
                 + \frac{2 T(w)}{ (z-w)^2 } 
                 + \frac{\partial T(w)}{z-w}    \nonumber \\
    T(z) W(w) & \sim & \frac{3W(w)}{(z-w)^2}
                     + \frac{\partial W(w)}{z-w} \label{ope_three}\\
    W(z) W(w) & \sim & \frac{-2/3}{(z-w)^6}
                     + \frac{2T(w)}{(z-w)^4}
                     + \frac{\partial T(w)}{(z-w)^3} \nonumber \\
                 &&  + \frac{1}{(z-w)^2}
                        \left(
                         \frac{8}{3}:T(w)T(w): - \hf {\partial}^2 T (w)
                        \right) \nonumber \\
                 &&  + \frac{1}{z-w} 
                        \left(
                           \frac{4}{3} \partial \left(
                                               :T(w) T(w):
                                               \right)
                           - \frac{1}{3} {\partial}^3 T (w)
                        \right).                    \nonumber 
  \end{eqnarray}

Representation theory of the vertex algebra $\VM$ is
just the same as that of $\/ {\cal U}(\Wth)$.
We see from the following lemma that $ \VMtwo $
is reducible so its maximal proper submodule $I$ is not zero. 
Representation theory of $\Wthtwo$ becomes
highly non-trivial due to the following constraints:
a module $M\/$ of the vertex algebra
$\VMtwo$ can be a module of
$\Wthtwo$ if and only if $M\/$ is annihilated
by all the Fourier 
components of all fields corresponding to vectors
in $I \subset \VMtwo$.

So it is important to find information of (top) singular vectors
in the vacuum module $\VMtwo$.
The following lemma can be proved by a tedious  
however direct calculation:
\begin{lemma}
 \begin{enumerate}
  \item[1)] There is no singular vector in $(\VMtwo)_n, \; n \leq 5.$
  
  \item[2)] There are two independent singular vectors
  in $(\VMtwo)_6$, denoted by $v_s$ and $v_s^{'}$:
   \begin{eqnarray*}
     v_s & \equiv& \left( \frac32 W_{-3}^2 
              - \frac{19}{36} L_{-3}^2 - \frac89 L_{-2}^3
              - \frac{14}{9} L_{-2}L_{-4} 
              + \frac{44}{9} L_{-6}  \right) \vac,          \\
     v_s^{'} & \equiv& \left( \frac92 W_{-6} + 9 L_{-3} W_{-3} -6 L_{-2} W_{-4}
               \right)  \vac.
    \end{eqnarray*}

  \item[3)] $v_s^{'} = \frac{27}{98} W_0 (v_s),
       \quad v_s = \frac{1}{36} W_0 (v_s^{'}).$ Equivalently we have

  \item[4)] $W_0 \left( 6 v_s \pm \frac{98}{27} v_s^{'} \right)
        = \pm 6 \left( 6 v_s \pm \frac{98}{27} v_s^{'} \right).$
  \end{enumerate}
  \label{lem_sing}
 \end{lemma}             

\begin{remark}
 Vectors $v_s, v_s^{'}$ are not singular vectors in the 
 Verma module $\Moo$.
\end{remark}

\section{Irreducible modules $V_{\alpha} (\alpha \in \Bbb C)$
of $\Wthtwo$}
  \label{sect_realize}

We first recall how we realize the $\Wthtwo$ algebra in terms
of a pair of fermionic $bc$ fields \cite{W2}.
Take a pair of $bc$ fields 
$$ b(z) = \SUM b (n) z^{-n}, \quad  c(z) = \SUM c (n) z^{-n-1} $$
with OPEs
\begin{equation}
b(z) c(w) \sim \frac{1}{z-w},
\quad b(z) b(w) \sim 0,
\quad c(z) c(w) \sim 0.
   \label{ope_bc}
\end{equation}
Equivalently, we have the following commutation relations:
$$ \left[ b(m), c(n) 
   \right]_{+} 
     = \delta_{m, -n}, \quad
   \left[ b(m), b(n)
   \right]_{+} = 0, \quad
   \left[ c(m), c(n)
   \right]_{+} =0.
$$
We denote by $\F$ the Fock space of the $bc$ fields, generated by
$\vacbc$, satisfying
$$ b( n+1 ) \vacbc = 0, 
   \quad  c(n) \vacbc = 0, \quad n \geq 0. $$
Then 
$$ j (z) = \/ : b(z) c(z) :\/ = \SUM j_n z^{-n-1} $$ 
is a free boson of conformal weight $1$ with commutation relations
$$[j_m, j_n] = m \delta_{m, -n}, \quad m,n \in \Z. $$ 
We further have the following commutation relations:
$$ [j_m, b(n)] = b(m+n), \quad
   [j_m, c(n)] = - c(m+n), \quad m,n \in \Z.
$$
Then we have the $bc$--charge decomposition of $\F$ 
according to the eigenvalues of $j_0$:
$$ \F = \bigoplus_{l \in \Z} \Fl .$$

We denote by $\Halpha$ ($\alpha \in \C$) the Fock space
of the Heisenberg algebra generated by $j_n, n \in \Z$,
with vacuum vector $\vacuum$ satisfying
$$ j_n \vacuum = \alpha \delta_{n,0}\vacuum, \quad n \geq 0.
$$
Denote by $\ps = q + j_0 \ln z - \sum_{n \neq 0} j_n z^{-n},$ 
where the operator $q$ satisfies
$[q, j_n] = \delta_{n,0}.$ Clearly $j (z) = \partial \ps.$
(Note that our $j(z), j_n, \cdots$ are denoted in \cite{W2}
by $-j^{bc} (z), -j^{bc}_n, \cdots$).

By the well-known boson-fermion correspondence, we have an 
isomorphism between $\Fl$ and $\Hl$ as representations over
the Heisenberg algebra generated by $j_n, n \in \Z$.
On the other hand, we may regard $b(z)$ and $c(z)$ as
\begin{equation}
  b(z) =: e^{\ps}:, \quad c(z) =: e^{- \ps}:.
   \label{eq_bf}
\end{equation}
Furthermore we have the following OPEs
\begin{equation}
  b(z) c(w) = \frac{1}{z-w} :b(z)c(w):, \quad 
  c(z) b(w) = \frac{1}{z-w} :c(z)b(w):.
   \label{ope_fermion}
\end{equation}
In particular it is well known that $\Fo$ (and so $\Ho$)
is a vertex algebra. 
Denote by $\Fbaro$ the kernel of the screening
operator $c(0)$ from $\Fo$ to ${\cal F}^{-1}$. It has
a structure of a vertex subalgebra of $\Fo$. Let 
\begin{equation}
  T (z) \equiv \SUM L_n z^{-n-2}
    = : \partial b (z) c(z) :. 
  \label{eq_vir}
\end{equation}
Easy to check that 
$T(z)$ is a Virasoro field with central charge $-2$.
We also define another field of conformal weight $3$:
\begin{equation}
  W(z) \equiv \SUM W_n z^{-n-3}
         = \frac{1}{\sqrt{6}}
              \left(
                : {\partial}^2 b(z) c(z) : 
                  - :\partial b(z) \partial c(z) :
              \right). 
   \label{eq_w3}
\end{equation}
We introduce a rescaled field 
$\widetilde{W} (z) = \hf \sqrt{6} \/ W(z)$
for convenience later on, namely
\begin{equation}
 \widetilde{W} (z)  \equiv \SUM 
  \widetilde{W}_n z^{-n-3}
   = \hf \left(
                : {\partial}^2 b(z) c(z) : 
                  - :\partial b(z) \partial c(z) :
         \right).
  \label{eq_wtilde}
\end{equation}

The following theorem is proved in \cite{W2}.
\begin{theorem}
  The vertex algebra $\Fbaro$ is isomorphic
 to the simple vertex algebra $\Wthtwo$ with
 generating fields $T(z)$ and $W(z)$ defined
 as in (\ref{eq_vir}) and (\ref{eq_w3}).
\end{theorem}

By the boson-fermion correspondence $\Fo$ and $\Ho$ are
isomorphic as vertex algebras so we may view $\Wthtwo$
as a vertex subalgebra of $\Ho$ as well. $\Halpha$
is a module over the vertex algebra $\Ho$ and so can 
be regarded as a module over $\Wthtwo$.
Denote by $V_{\alpha}$ the irreducible subquotient
of the $\Wthtwo$-submodule of $\Halpha$ generated
by the highest weight $\vacuum$. 

We first rewrite the fields $T(z)$ and $W(z)$ defined
 as in (\ref{eq_vir}) and (\ref{eq_w3}) in terms of
the field $j(z)$ and its derivative fields.
\begin{proposition}
  Under the boson-fermion correspondence, the
fields $T(z)$ and $W(z)$ in (\ref{eq_vir}) and (\ref{eq_w3}) 
can be expressed in terms of $j(z)$ as
  \begin{eqnarray}
   T(z) &=& \hf \left( :j(z)^2: + \partial j(z) \right),  
                                       \label{eq_tboson}  \\
   \widetilde{W} (z) &=& \frac{1}{12} 
    \left( 4 :j(z)^3 : + 6 :j(z) \partial j(z): + \partial^2 j(z) 
    \right).
    \label{eq_wboson}
  \end{eqnarray}
  \label{prop_bosonw3}
\end{proposition}
\begin{demo}{Proof}
 By (\ref{ope_fermion}), we have
 \begin{eqnarray*}
  b(z) c(w) & \sim &
   \zwminusone
    \left\{
     1 + \zwone j(w) + \frac{\zwtwo}{2} \left(
                             :j(w)^2: + \partial j(w) 
                                  \right) \right.     \\
       && 
    \left.
      + \frac{\zwthree}{6} \left( :j(w)^3: + 3 :j(w) \partial j(w): + 
                            \partial^2 j(w) 
                         \right)
    \right\}                    \\
            & \sim &
     \zwminusone + j(w) + \hf \zwone \left( 
                         :j(w)^2: + \partial j(w) 
                                     \right)            \\
      & & + \frac16 \zwtwo \left( 
       :j(w)^3: + 3 :j(w) \partial j(w): + \partial^2 j(w) 
                           \right).
 \end{eqnarray*}
From this we see that 
\begin{equation}
 : \partial^2 b(w) c(w) :
    = \frac13 \left(:j(w)^3: + 3 :j(w) \partial j(w): + 
                            \partial^2 j(w)
              \right).
  \label{eq_normal}
\end{equation}
On the other hand, by (\ref{ope_fermion}) we have 
 \begin{eqnarray*}
 c(z) b(w) &= &
  \zwminusone 
   \left\{
     1 - \zwone j(w) + \frac{\zwtwo}{2} \left(
                             :j(w)^2: - \partial j(w) 
                                  \right)      
   \right.     \\
     & & 
   \left.  + \frac{\zwthree}{6} \left(
       - :j(w)^3: + 3 :j(w) \partial j(w): - \partial^2 j(w) 
                              \right)
   \right\}   \\
   & & + \mbox{ higher terms}.     
  \end{eqnarray*}
This implies that 
 \begin{eqnarray}
  :\partial c(z) b(w): & = &
   - \zwminustwo + \hf (:j(w)^2: - \partial j(w) )  \nonumber      \\
   & & + \frac13 \zwone \left( 
             - :j(w)^3: + 3 :j(w) \partial j(w): 
                        \right)    \nonumber \\
   & & + \mbox{higher terms}.
  \label{eq_dcb}
 \end{eqnarray}

Equivalently we have by switching $ z$ and $w$ in (\ref{eq_dcb})
and reversing the order between $\partial c$ and $b$ (we get a 
minus sign since $bc$ fields are fermionic)
 \begin{eqnarray}
  :b(z) \partial c(w):
   & \sim &
     \zwminustwo - \hf \left( :j(z)^2: - \partial j(z) 
                       \right)       \nonumber         \\
     & &   + \frac13 \zwone \left(   
                       - :j(z)^3: + 3 :j(z) \partial j(z): 
                            \right)   \nonumber \\
   & \sim &
     \zwminustwo - \hf \left( 
                        :j(w)^2: - \partial j(w) 
                       \right)  \nonumber \\    
     &&  - \hf (z-w) \left( 
                      2: j(w) \partial j(w): - \partial^2 j(w) 
                  \right)        \nonumber \\
     && + \frac13 \left( 
                    - :j(w)^3: + 3 :j(w) \partial j(w): 
                  \right)   \nonumber \\
   & \sim &
     \zwminustwo - \hf \left( :j(w)^2: - \partial j(w) 
                       \right)   \nonumber \\
      && - \hf (z-w) \left( 2: j(w) \partial j(w): - \partial^2 j(w)
                     \right )    \nonumber \\
      && + \zwone \left(
                 -\frac13 :j(w)^3: + \frac16 \partial^2 j(w)      \right).
    \label{ope_dbdc}
  \end{eqnarray}
It follows from (\ref{ope_dbdc}) that
 \begin{equation}
  : \partial b(w) \partial c(w):
    = -\frac13 :j(w)^3: + \frac16 \partial^2 j(w).
   \label{eq_nop}
 \end{equation}
So by (\ref{eq_wtilde}), (\ref{eq_normal}) and (\ref{eq_nop}) we have 
$$ \widetilde{W} (w) =
  \frac{1}{12} 
    \left( 4 :j(w)^3 : + 6 :j(w) \partial j(w) + \partial^2 j(w) 
    \right).
$$
 The proof of the identity (\ref{eq_tboson}) is similar.
\end{demo}

\begin{proposition}
  The highest weight of the $\Wthtwo$-module
$V_{\alpha}$ ($\alpha \in \C$)
with respect to $(L_0, \widetilde{W}_0)$ is
$\left( \hf \alpha (\alpha -1),
         \frac16 \alpha (\alpha -1)( 2 \alpha -1 )
 \right).$
\end{proposition}

\begin{demo}{Proof}
   $L_0$ and $\widetilde{W}_0$ can be written as an infinite 
sum of monomials in terms of $j_{-n}, n > 0$ by Proposition
\ref{prop_bosonw3}. Indeed we have
\begin{eqnarray*}
  L_0 &=& \hf \left(
              \sum_{n < 0} j_n j_{-n} + \sum_{n \geq 0} j_{-n} j_n
            \right)  -\hf j_0    \\
      &=& \hf (j_0^2 - j_0) + \sum_{n > 0} j_{-n} j_n.
\end{eqnarray*}
Since $j_n \vacuum = 0, n > 0 $ and $j_0 \vacuum = \alpha \vacuum$,
we have 
$$L_0 \vacuum = \hf (j_0^2 - j_0) \vacuum 
=\hf \alpha (\alpha -1) \vacuum. $$
Similarly, a little calculation shows that
the only terms in $\widetilde{W}_0$ which do not 
annihilate the vacuum vector $\vacuum$ are 
$ \frac{1}{12} \left( 4 j_0^3 + 6 j_0 ( - j_0 ) + 2 j_0 \right). $
So we have 
\begin{eqnarray*}
\widetilde{W}_0 \vacuum 
   = \frac{1}{12} \left(
          4 {\alpha}^3 + 6 \alpha (- \alpha) + 2\alpha 
                  \right)
   = \frac16 \alpha (\alpha -1)( 2 \alpha -1 ).
\end{eqnarray*}
\end{demo}

\begin{remark}
  The irreducible module $V_{\alpha}$ is isomorphic to the module
$\overline{\cal F}^{- \alpha}$ constructed 
in \cite{W2} by comparing their highest weights for $\alpha \in \Z$.
$V_{\alpha}$ is a proper subspace of $\Halpha$ in this case
and its full character formula is given in \cite{W2}.
For $\alpha \not \in \hf \Z$,
we know that $\Halpha$ is irreducible as a module
over the Virasoro algebra given by the field $T(z)$ with central
charge $-2$ \cite{FeF, KR} and so is irreducible as a module
over $\Wthtwo$. Full character formulas of these $V_{\alpha}$
with respect to $\{ L_0, \widetilde{W}_0 \}$ can be 
also calculated.
\end{remark}

\section{Classification of irreducible representations of $\Wthtwo$ algebra}
  \label{sect_complete}

We will show that irreducible modules $V_{\alpha}, \alpha \in \C$
exhaust all the irreducible modules of vertex algebra $\Wthtwo$ by
calculating Zhu algebra in this case. We break the proof into
a sequence of simple lemmas.
\begin{lemma}
   Zhu algebra $A(\VM)$ is isomorphic to a polynomial
  algebra $\C [t, w]$, where $t, w$ correspond to
  $[L_{-2} \vac ]$ and $[\widetilde{W}_{-3} \vac ]$
  in $A(\VM )$. 
   \label{lem_zhuverma}
\end{lemma}             
Note that $ L_{-2}\vac$ is the Virasoro element in $\VM$
so the element $[L_{-2}\vac ]$ lies in the center 
of $A(\VM)$ by Lemma \ref{lem_center}.
Proof of the above lemma is quite standard. See Lemma 4.1 in
\cite{W1} for a proof of a similar result. One can easily modify
that proof to give a proof of our present lemma. We will not
write it down here since it is not very illuminating.

Now specify $c = -2$.
Let us denote by $\sigma$ the isomorphism 
from $A(\VMtwo)$ to $\C [t, w]$. 
\begin{lemma}
Keeping the conventions in Lemma \ref{lem_zhuverma},
under the isomorphism $\sigma$ we have
  $$\sigma {( {[ v_s ]} )} = w^2 - \frac19 t^2 (8t + 1),\quad 
    \sigma {( {[ v^{'}_s ]} )} = 0.$$
   \label{lem_zhusing}
\end{lemma}
\begin{demo}{Proof}
  We will continue using the equivalence convention 
  denoted by $\sim$ in the sense of Section \ref{sect_general}.
  It follows from lemma \ref{lem_keylemma} that for any $a \in \VMtwo$,
  \begin{eqnarray}
   a * \left( \widetilde{W}_{-3} \vac
             \right) 
   & \sim & \left(
           \widetilde{W}_{-3} + 2 \widetilde{W}_{-2} 
           + \widetilde{W}_{-1} 
         \right) a  \nonumber \\
     a * \left( L_{-2} \vac
         \right) 
   & \sim & \left(
           L_{-2} + L_{-1}  
         \right) a.     
    \label{eq_equi}
  \end{eqnarray}
  Recall that the isomorphism $\sigma$ from $A(\VMtwo)$ to
  $\C [t, w]$ sends elements $[L_{-2} \vac ]$ and
  $[\widetilde{W}_{-3} \vac ]$
  in $A(\VM )$ to $t, w$ respectively.
  By applying (\ref{eq_equi}) to the first two terms of
  the singular vector $v_s$ given in Lemma \ref{lem_sing} and then rewriting
  it in terms of the PBW basis of the form (\ref{eq_PBW})
  by using the commutation relations (\ref{eq_alg}),
  we get 
  \begin{eqnarray}
    v_s  & \sim & w^2
           + \left\{ - \left(
                             2\widetilde{W}_{-2} 
                              + \widetilde{W}_{-1} 
                       \right)\widetilde{W}_{-3}
             \right.      \nonumber        \\
          &  & \left.  
                     - \frac{19}{36} L_{-3}^2 - \frac89 L_{-2}^3
                     - \frac{14}{9} L_{-2}L_{-4} 
                     + \frac{44}{9} L_{-6}  
                   \right\} \vac \nonumber        \\
         & = & w^2
           + \left\{ - 3 \left(
                     \frac83 L_{-2}L_{-3} - \frac{10}{3} L_{-5}
                         \right) \right.       \nonumber           \\
         &  &  - \frac32 \left( \frac83 L_{-2}^2 - L_{-4} 
                         \right)  \nonumber  \\
         &  &  \left. - \frac{19}{36} L_{-3}^2 - \frac89 L_{-2}^3
              -\frac{14}{9} L_{-2}L_{-4}
             + \frac{44}{9} L_{-6}
           \right\} \vac      \nonumber  \\
         & = & w^2 +\left\{ 
                     - 8 L_{-2}L_{-3} + 10 L_{-5} 
                     -4 L_{-2}^2 + \frac32 L_{-4}  \right. \nonumber  \\
         && \left. - \frac{19}{36} L_{-3}^2 - \frac89 L_{-2}^3
              -\frac{14}{9} L_{-2}L_{-4}
             + \frac{44}{9} L_{-6}  \right\} \vac .
     \label{eq_iden}  
  \end{eqnarray}
  It is easy to show by induction and applying Lemma \ref{lem_keylemma}
  that
  \begin{equation}
     L_{-n} \sim (-1)^n \left(
                      (n-1)\left( L_{-2} + L_{-1} \right) + L_0 \right),
     \quad n \geq 1.
   \label{eq_co}
  \end{equation}
  By the equation (\ref{eq_equi}) and repeated uses of 
 (\ref{eq_co}) on the right hand side of (\ref{eq_iden}), we get
  \begin{eqnarray}
     v_s & \sim & w^2 + 16 t(t+3) - 40t -4t(t+2) 
        + \frac92 t  \nonumber          \\
       & & - \frac{19}{18} t(2t+3) - \frac89 t(t+2)(t+4) 
            - \frac{14}3 t(t+4) +\frac{220}{9} t    \nonumber          \\
       &= & w^2 - \frac19 t^2 (8t + 1).
  \end{eqnarray}

  This completes the proof that 
  $\sigma {( {[ v_s ]} )} = w^2 - \frac19 t^2 (8t + 1).$
  Similarly we can prove that
  $ \sigma {( {[ v^{'}_s ]} )} = 0.$
\end{demo}

Denote $ f(t, w) = w^2 -\frac19 t^2 (8t + 1).$
Now the following lemma follows from 
Proposition \ref{prop_associativity},
Lemma \ref{lem_zhuverma} and Lemma \ref{lem_zhusing} 
(see Corollary \ref{cor_exact}
for a more precise statement).
\begin{lemma}
  The Zhu algebra $A(\Wthtwo )$ is a certain quotient of the quotient
  algebra $\C [t, w] / <f(t, w) >$,
  where $<f(t, w) >$ denotes the ideal of $\C [t, w]$
  generated by $f(t,w) \in \C [t, w]$.
 \label{lem_key}
\end{lemma}             
We have the following observation.
\begin{lemma}
  Solutions to the equation (\ref{eq_elliptic}) are parametrized 
  as follows:
  \begin{equation}
   \left(  t(\alpha), w(\alpha) \right) \equiv
   \left( \hf \alpha(\alpha -1), \frac16 \alpha(\alpha -1)(2\alpha -1)
   \right), \quad \alpha \in \C.
    \label{eq_solu}
  \end{equation}
   \label{lem_solu}
\end{lemma} 
\begin{demo}{Proof}
  First it is clear that $t(\alpha)$ can take any complex value 
 when $\alpha$ ranges over $\C$. Then by substituting $t(\alpha)$
in the equation (\ref{eq_elliptic}) we see that
$w(\alpha)^2 = \left[ \frac16 \alpha(\alpha -1)(2\alpha -1) \right]^2.$
We don't lose any generality by letting 
$w(\alpha) = \frac16 \alpha(\alpha -1)(2\alpha -1)$. The
reason is that $t(1 - \alpha) = t(\alpha) $ while 
$w(1- \alpha) = -w(\alpha) $.
\end{demo} 
    
\begin{remark}
  For different $\alpha, \alpha' \in \C$, 
  $\left(  t(\alpha), w(\alpha) \right) =
  \left(  t(\alpha'), w(\alpha') \right)$ if and only if 
  $\alpha = 0 (\mbox{resp. } 1),
  \alpha' = 1(\mbox{resp. } 0).$ Namely $V_0$ is isomorphic
  to $V_1$ and this is the only isomorphism 
  among $V_{\alpha}, \alpha \in \C.$
 \label{rem_iso}
\end{remark}
 
Now we are ready to prove our classification theorem on irreducible
modules over the $\Wthtwo$ algebra.
\begin{theorem}
  $V_{\alpha}, \alpha \in \C$ are all the irreducible
  modules over the simple $\Wth$ algebra with central charge
  $-2$. Highest weights of these modules $V_{\alpha}$
  are given by $\left( \hf \alpha(\alpha -1), 
  \frac16 \alpha(\alpha -1)(2\alpha -1)
   \right), \alpha \in \C.  $
  They are parametrized by points $ (t, w)$ on the rational curve 
  defined by $w^2 = \frac19 t^2(8t + 1).$
   \label{th_main}
\end{theorem}
\begin{demo}{Proof}
   By Lemma \ref{lem_key}, we see that any irreducible
module of the associative algebra $A(\Wthtwo)$
is one-dimensional since $A(\Wthtwo)$ is commutative. 
Given $t, w \in \C$, let $\C_{t,w}$ 
be the one-dimensional module
of $A(\Wthtwo)$, with $[L_{-2} \vac]$ acting as the scalar
$t$ and $[\widetilde{W}_{-3} \vac]$ as the scalar
$w$. Then $(t, w)$ has to satisfy 
$w^2 = \frac19 t^2(8t + 1).$ 
Note that $o(L_{-2} \vac ) = L_0$ and 
$o( \widetilde{W}_{-3} \vac ) = \widetilde{W}_0$ by the 
definition of $o(\cdot )$ in Section~\ref{sect_general}.
So by Theorem~\ref{th_zhucorresp},
the highest weight $(t, w)$ of any irreducible module
of the vertex algebra $\Wthtwo$
with respect to $(L_0, \widetilde{W}_0)$ has to satisfy
the equation $w^2 = \frac19 t^2(8t + 1).$
By Lemma \ref{lem_solu}, we see all solutions to the above
equation can be written as of the form
$\left( \hf \alpha(\alpha -1), \frac16 \alpha(\alpha -1)(2\alpha -1)
   \right),$ $ \alpha \in \C.  $
On the other hand, we have already constructed a family of 
irreducible modules $V_{\alpha}\;( \alpha \in \C )$ with
highest weight exactly equal to
$\left( \hf \alpha(\alpha -1), \frac16 \alpha(\alpha -1)(2\alpha -1)
   \right).  $
This completes the proof of the theorem.
\end{demo}

We think it is remarkable that the set of all irreducible
modules of $\Wthtwo$ has such a simple and elegant description
in terms of a rational curve. It indicates that the non-rational
vertex algebras may have very rich representation theory.

We have an immediate corollary of Theorem \ref{th_main} which
strengthens Lemma \ref{lem_key}. 

\begin{corollary}
  The Zhu algebra $A(\Wthtwo)$ is isomorphic to the
commutative associative algebra $\C [t, w] / <f(t, w)>$.
  \label{cor_exact}
\end{corollary}

\begin{remark}
\begin{enumerate}
  \item[1)] Basing on the results of
  Theorem \ref{th_main} and Corollary~\ref{cor_exact}
  it is natural to conjecture that the singular vectors $v_s$ and
  $v^{'}_s$ generate the maximal proper submodule of
  the vacuum module $\VMtwo$.

  \item[2)] A Virasoro vertex algebra with a certain central charge 
is rational
if and only if the corresponding vacuum
module is reducible \cite{W1}. As our results show, ${\cal W}_3$
algebra provides new possibility, namely the simple vertex algebra
$\Wthtwo$ is not rational but the corresponding vacuum module
$\VMtwo$ is reducible.
  \end{enumerate}
\end{remark}

We further comment on why central charge $c= -2$ is particularly
interesting from a different point of view. 
There is the so-called quantized Drinfeld-Sokolov
reduction (cf. e.g. \cite{BH, FKW})
which allows one to establish connections
between ${\cal W}_n$ algebra with central charge $c_n^{(k)}$ and 
the affine Kac-Moody Lie algebra $\widehat{sl}_n$ with central
charge $k$. Here 
$$ c_n^{(k)} = 2n^3 -n -1 -n(n^2 -1)
    \left(
        \frac{1}{k+n} + k +n 
    \right).$$

In particular, for $k = -n + p/q$ one can rewrite $ c_n^{(k)}$
as follows:
$$ c_n^{(k)} = (n-1) \left(
                     1- \frac{n(n+1)(p-q)^2}{pq}
                  \right).$$
The so-called minimal series central charges of the ${\cal W}_n$ algebra
are those
$ c_n^{(k)}$ for $k = -n + p/q,$
where $p, q$ are coprime integers satisfying 
$ p, q \geq n.$ They correspond to 
the admissible central charges $k = -n + p/q$
for the affine algebra $\widehat{sl}_n$,
with the same conditions imposed on $p, q$ as above.
The admissible representations with admissible central
charges were first studied by Kac-Wakimoto \cite{KWa}.

Thus by means of Drinfeld-Sokolov reduction
the central charge $-2$ for the $\Wth$ algebra
corresponds to the central charge $k = -\frac32$ or $-\frac72$
of $\widehat{sl}_3$.
Observe that $k = -\frac32 = -3 + \frac23$ or $-\frac72 = -3 + \frac32$ 
corresponds to the
``boundry'' of the admissible central charges of $\widehat{sl}_3$.

However more than this is true. Consider the                            
``boundry'' of the admissible central charges of $\widehat{sl}_n$,
i.e. $ k = -n + \frac{n}{n-1}$ or $ -n + \frac{n-1}{n}$.
The corresponding central charge of the ${\cal W}_n$ algebra
$c_n^{(k)} = -2$, which is independent of $n$. 
In this sense $-2$ is a universal 
central charge for any ${\cal W}_n$ algebra. We expect that the
representations of the affine algebra $\widehat{sl_n}$
with central charge equal to the ``boundry'' of the admissible
central charges are of independent interest.

\section{Classification of irreducible modules of
vertex algebra $\Winfone$}
\label{sect_infty}
  
Let $\D$ be the Lie algebra of regular differential operators on
the circle. The elements 
\begin{eqnarray*}
  J^l_k = - t^{l+k} ( \partial_t )^l, 
     \quad l \in \Z_{+}, k \in \Z, 
\end{eqnarray*}
form a basis of $\D$. $\D$ has also another basis
\begin{eqnarray*}
  L^l_k = - t^{k} D^l, 
     \quad l \in \Z_{+}, k \in \Z,  
\end{eqnarray*}
where $D = t \partial_t$. Denote by $\hD$ the central extension of  
$ \D $
by a one-dimensional center with a generator $C$, with
commutation relation (cf. \cite{KR1})
\begin{eqnarray}
  \left[
     t^r f(D), t^s g(D)
  \right]
    & = & t^{r+s} 
    \left(
      f(D + s) g(D) - f(D) g(D+r) 
    \right) \nonumber \\ 
   & + & \Psi 
      \left(
        t^r f(D), t^s g(D)
      \right)
      C,
  \label{eq_12}
\end{eqnarray}
where 
\begin{equation}
  \Psi 
      \left(
        t^r f(D), t^s g(D)
      \right)
   = 
   \left\{
      \everymath{\displaystyle}
      \begin{array}{ll}
        \sum_{-r \leq j \leq -1} f(j) g(j+r),& r= -s \geq 0  \\
        0, & r + s \neq 0. 
      \end{array}
    \right. \\
  \label{eq_13}
\end{equation}

Letting weight $J^l_k = k$ and weight $ C = 0$ defines a principal
gradation
\begin{equation}
  \hD = \bigoplus_{j \in \Z} \widehat{\cal D}_j.
  \label{eq_14}
\end{equation}
Then we have the triangular decomposition
\begin{equation}
  \hD = \hD_{+} \bigoplus \hD_{0}  \bigoplus \hD_{-},
  \label{eq_15}
\end{equation}
where 
\begin{eqnarray*}
  \hD_{\pm} = \bigoplus_{j \in \pm \Bbb N} \widehat{\cal D}_j,
  \quad
  \hD_{0} = \D_{0}  \bigoplus {\Bbb C} C.
\end{eqnarray*}

Let $\cal P$ be the distinguished parabolic subalgebra of $\D$,
consisting of the differential operators that extends
into the whole interior of the circle. 
${\cal P}$ has a basis $\{ J^l_k, l \geq 0, l + k \geq 0 \}.$
It is easy to check that the $2$-cocycle $\Psi$ defining
the central extension of $\hD$ vanishes when restricted
to the parabolic subalgebra $\cal P$. So $\cal P$ is also a
subalgebra of $\hD$. Denote
$\widehat{\cal P} = {\cal P} \oplus {\Bbb C} C$.

Fix $c \in \Bbb C$. Denote by ${\Bbb C}_c$ the $1$--dimensional
$\widehat{\cal P}$ module by letting $C$ acts as scalar $c$ and 
$\cal P$ acts trivially. Fix a non-zero vector
$v_0$ in ${\Bbb C}_c$. The induced $\hD$--module
\begin{eqnarray*}
  M_c \left( \hD
      \right)
  = {\cal U} \left( \hD
             \right)
    \bigotimes_{ {\cal U} \left( \cal P
                          \right)
               }
    {\Bbb C}_c
\end{eqnarray*}
is called the vacuum $\hD$--module with central charge $c$. Here we 
denote by ${\cal U} (\frak g)$ the universal enveloping
algebra of a Lie algebra $\frak g$. $M_c (\hD)$ admits a unique
irreducible quotient, denoted by $\Wc$. 
Denote the highest weight vector $1\otimes v_0$ in
$M_c (\hD)$ by $\vac$.

It is shown in 
\cite{FKRW} that $\Wc$ carries a canonical vertex
algebra structure, with vacuum vector $\vac$ and
generating fields 
\begin{eqnarray*}
  J^l (z) = \sum_{k \in \Bbb Z} J^l_k z^{-k-l-1},
\end{eqnarray*}
of conformal weight $l + 1, l = 0, 1, \cdots.$ The fields
$J^l (z)$ corresponds to the vector $J^l_{-l-1} \vac$
in $\Wc$. Below we will concentrate on the 
particular case $\Winfone$.
       
The relation between vertex algebras $\Winfone$ and
$\Wthtwo$ is made clear by the following theorem \cite{W2}.

\begin{theorem}
  The vertex algebra $\Winfone$ is isomorphic to a tensor
product of the $\Wthtwo$ algebra, and the Heisenberg 
vertex algebra ${\cal H}_0$ with $ J^0 (z)$ as 
a generating field.
\end{theorem}

Then the classification of irreducible modules over
$\Winfone$ follows from classification of those
over $\Wthtwo$ since the classification of irreducible
modules over a Heisenberg vertex algebra is well known.
Also see Remark \ref{rem_iso}.

\begin{theorem}
  There exists a two-parameter family of irreducible modules
over $\Winfone$. Any irreducible $\Winfone$-module
can be written uniquely as a tensor product of a module
${\cal L}(t(\alpha), w(\alpha))$ 
over $\Wthtwo$ with a module ${\cal H}_s$ over ${\cal H}_0$
($ \alpha\in \C-\{1\}$, $s \in \C$), with 
$(t(\alpha), w(\alpha))$ as defined in (\ref{eq_solu}).
  \label{th_two}
\end{theorem}

\begin{remark}
   Theorem \ref{th_two} disproves a conjecture of Kac and Radul
  \cite{KR2}. The list of irreducible modules of $\Winfone$
  which were conjectured to be complete in \cite{KR2}
  consists of those with $\alpha = 0$
  in Theorem~\ref{th_two}, (i.e. modules ${\cal M}_s^0$ in \cite{W2}).
\end{remark}

There are several questions which the author does not know the
answers at present but hope to have a better understanding
in the near future.

\begin{enumerate}
   \item[1)] What are the fusion rules of 
  $\Wthtwo$ (and thus of $\Winfone$)? 
  The existence of reducible however indecomposible
  modules of $\Wthtwo$ \cite{W2} seems to be related to
  the fact that there is a node at $(0, 0)$ on  the
  rational curve $w^2 - \frac{1}{9} t^2 (8t + 1) = 0.$
  It is likely that we may need to regard some reducible
  however indecomposible modules as basic objects when
  studying the fusion rules.

  \item[2)] Recall the Cartan subalgebra of $\W$ is infinite
  dimensional. In \cite{KR1} the highest weight of an irreducible quasifinite
 module over $\W$ is characterized in terms of a certain generating
 function $\Delta (x)$. The question is how to identify
 highest weights of the two-parameter family of irreducible
 modules of $\W$ with central charge $-1$ in Theorem 
 \ref{th_two} in terms of $\Delta (x)$. It would be very interesting
 to see if these irreducible modules of $\Winfone$ we have constructed
 are the first realizations of $\W$-modules with 
 $\Delta (x) = \frac{ p(x) e^{sx}}{e^x -1} + \cdots$ 
 with some non-constant polynomial $p(x)$ (cf. \cite{KR1} for notations).
\end{enumerate}

\noindent{\bf Acknowledgement} Some results of this paper
were presented in the Seminar of
Geometry, Symmetry and Physics at Yale university 
and in the 1997 AMS Meeting at Detroit.
I thank the organizers of the meeting, C. Dong and R. Griess,
for invitation. I also thank E. Frenkel, I. Frenkel,
V. Kac and G. Zuckerman for their interests and comments.

\end{document}